\documentclass{article}
\usepackage[utf8]{inputenc}
\usepackage{graphicx}
\usepackage{url}
\usepackage{amsfonts}
\usepackage{hyperref}
\usepackage{eurosym}
\begin{document}

\title{Ellipsograph of Archimedes\\ as a simple LEGO construction}
\author{Sara Jodlbauer and Zoltán Kovács}
\date{Linz School of Education\\
Johannes Kepler University, Linz, Austria\\
A-4040 Linz, Altenberger Straße 69}

\maketitle

\begin{abstract}
    We report on a simple LEGO construction that can draw an ellipse by using the concept of trammel of Archimedes.
\end{abstract}

School curriculum related to elementary  geometry usually focuses on ``compass and ruler'' constructions where the traditional Euclidean construction steps are allowed to generate a figure. While it is possible to construct any point of an ellipse that is given by its parameters, school curriculum seldom covers any methods how to efficiently draw an arc or the full ellipse by home-made mechanical constructions.

To draw an ellipse on your own, several options exist. One of them, the gardener method is mentioned in many textbooks (see \cite{wiki:ellipse} for a reference), and it is based on the locus definition, namely:
\begin{quote}
Given two fixed points $F_1$, $F_2$  called the foci and a distance $2a$ which is greater than the distance between the foci, the ellipse is the set of points $P$ such that the sum of the distances $|PF_1|$, $|PF_2|$  is equal to $2a$:
$$E=\{P\in \mathbb {R} ^{2}:\,|PF_{1}|+|PF_{2}|=2a\}.$$
\end{quote}
By using this definition, the gardener puts down two stakes and loops a piece of rope around them. Using a stick, he pulls the loop taut and marks the points around a curve.

In this article we consider a different method, namely the ellipsograph (or trammel) of Archimedes \cite{Artobolevskii} that consists of two shuttles which are confined to perpendicular channels or rails and a rod which is attached to the shuttles by pivots at fixed positions along the rod.

We remark that such constructions are well-known and very popular as commercial toys. Also as LEGO constructions there are several options. A recent YouTube video at \url{https://www.youtube.com/watch?v=E5V9C-9dxUo} contributed by \textit{Brixe63} shows a complicated mechanism that is definitely based on the same idea, but consists of several LEGO bricks. Another approach (among several others) can be found on Rebrickable at \url{https://rebrickable.com/mocs/MOC-4096/JKBrickworks/trammel-of-archimedes/#details}, contributed by \textit{JKBrickWorks}, it uses 47 parts. In our contribution---which is based on \cite{jodlbauer-ba}---we focus on minimizing the number of parts, and, at the same time, we focus on \emph{drawing} the ellipse, not just the building up the motion.

\begin{figure}[!htbp]
\begin{center}
    \includegraphics[width=0.5\textwidth]{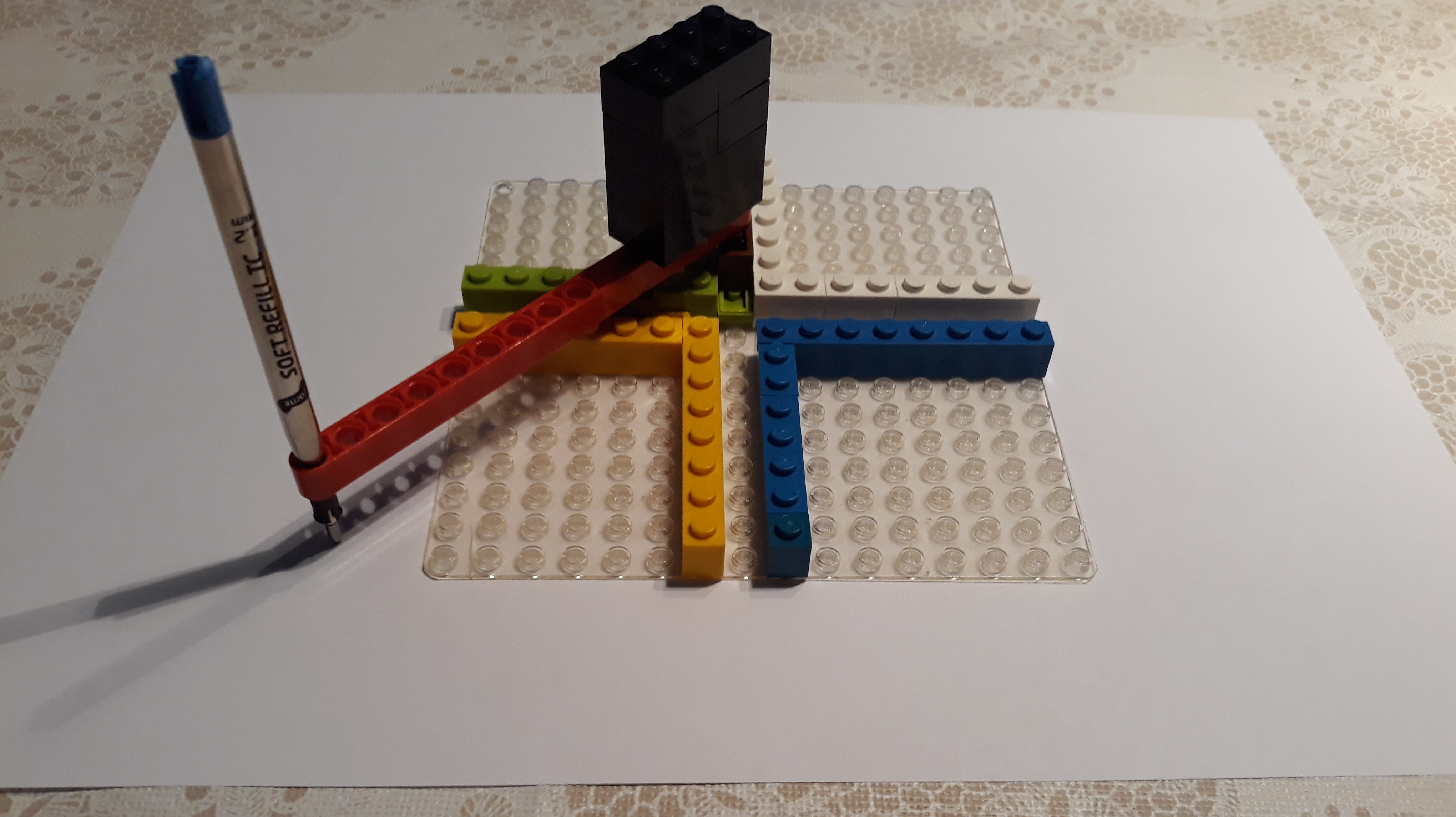}
    \caption{An ellipsograph}
    \label{ellipsograph1}
\end{center}
\end{figure}

\begin{figure}[!htbp]
\begin{center}
    \includegraphics[width=0.5\textwidth]{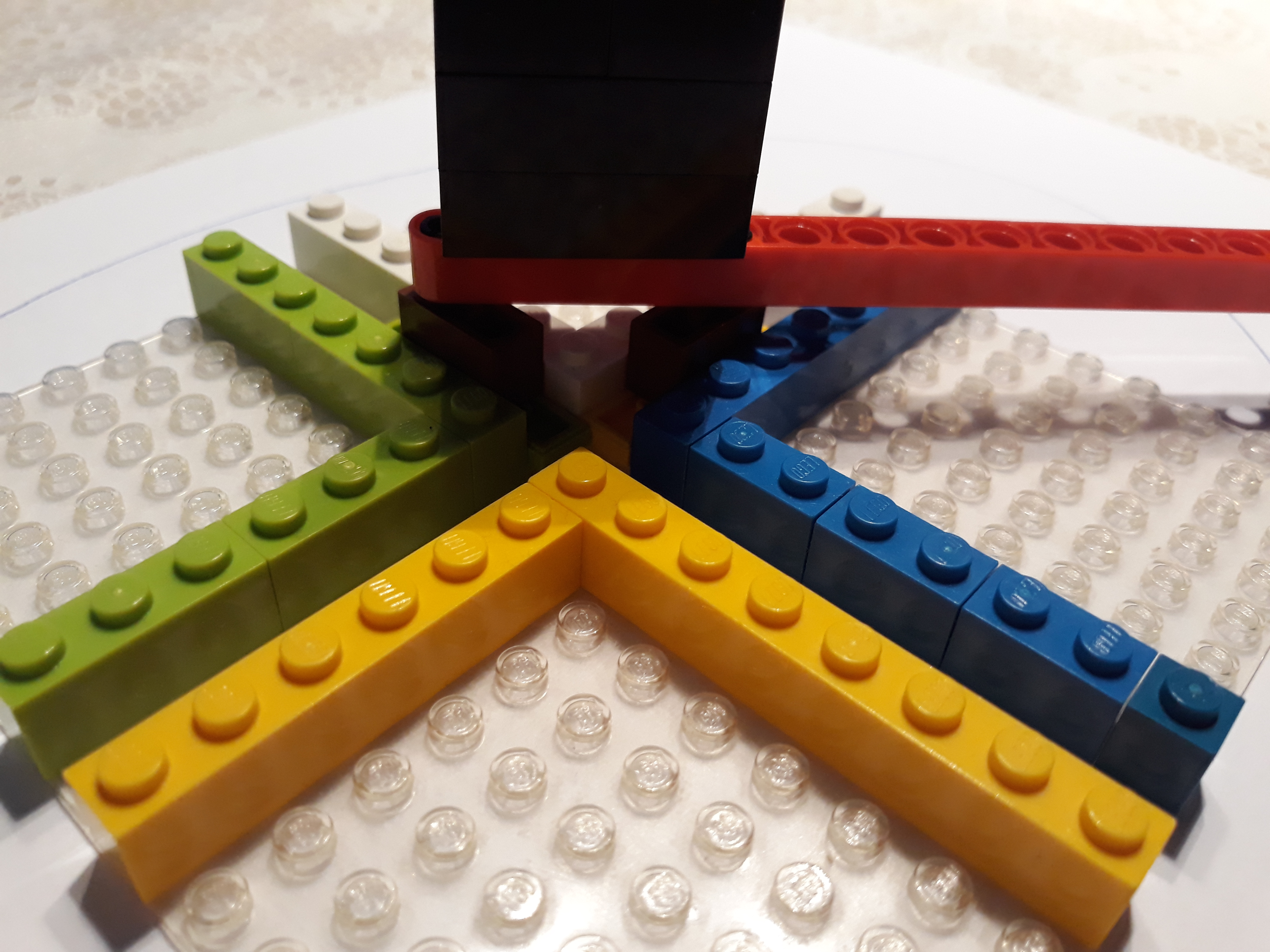}
    \caption{A closer look on the ellipsograph}
    \label{ellipsograph2}
\end{center}
\end{figure}

Fig.~\ref{ellipsograph1} shows how the ellipsograph is built. Easy sliding is assured by using two flat tiles for the shuttles. Their collision-free arrival in the channels is solved by using  $1\times4$ flat tiles. However, their somewhat long size has some drawback, namely, that the shuttles require a bigger distance to each other, as seen in Fig.~\ref{ellipsograph2}.

\begin{figure}[!htbp]
\begin{center}
    \includegraphics[width=1\textwidth]{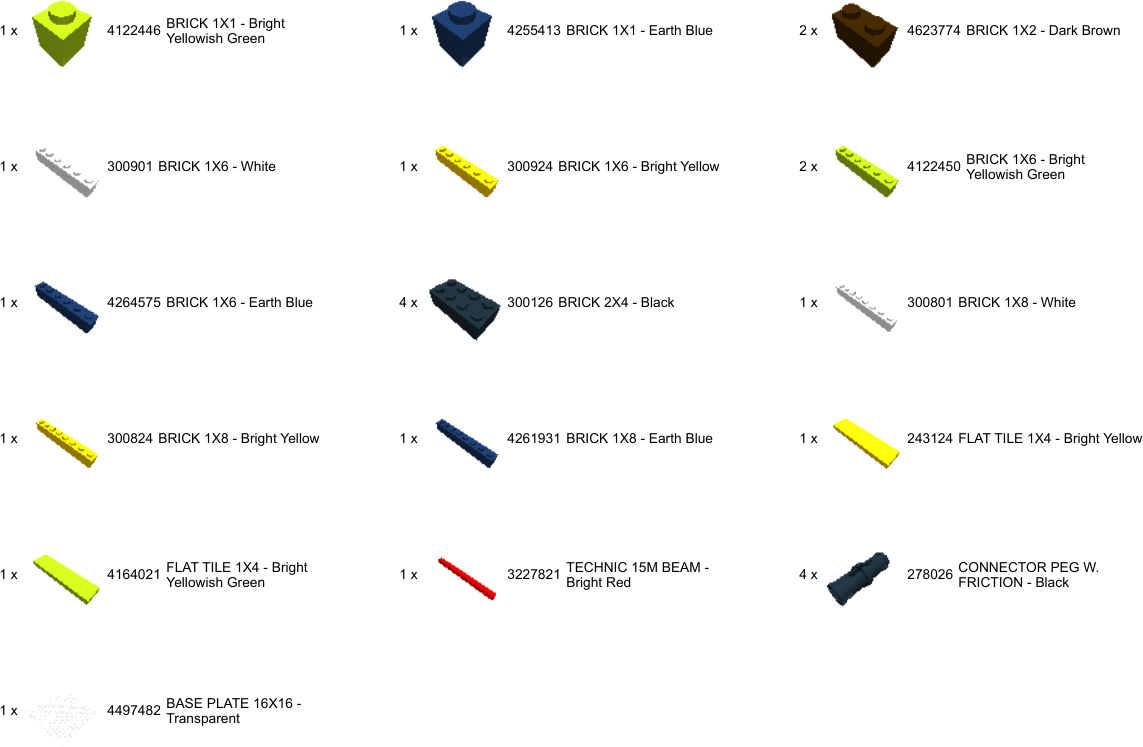}
    \caption{LEGO bricks used for the ellipsograph}
    \label{ellipsograph3}
\end{center}
\end{figure}

We use 24 LEGO parts listed in Fig.~\ref{ellipsograph3} to build the construction and a G2 type pen refill (see \url{https://en.wikipedia.org/wiki/Ballpoint_pen#/media/File:Ball_point_pen_refills_en.png} for an overview on pen refills). We highlight that this construction well harmonizes with a set of other LEGO bricks, described in \cite{lego} at the GitHub page \cite{gh-lego-linkages}; it allows the students to draw a high variety of other algebraic curves. On the other hand, static balancing of our contribution is an important part: a counterweight is required to ensure stable motion (see Fig.~\ref{ellipsograph3}).

We emphasize that the LEGO ellipsograph can be cheaply built by ordering parts from Internet stores. Table \ref{tab:set} shows some recent prices at \url{brickowl.com}, by selecting the cheapest color variant of the required parts at the time of writing this paper. Teachers and students should be aware that a suitable pen refill is also required, but it should be available at an affordable price below 1\euro . Thus the total minimal price of our ellipsograph is about 2\euro.

\begin{table}
\begin{center}
{\scriptsize
\begin{tabular}{ccccccccc|c}
     &     &     &     &     &Flat&   &   &Base&\\
Brick&Brick&Brick&Brick&Brick&tile&   &Pin&plate&\\
$1\times1$&$1\times2$&$1\times6$&$2\times4$&$1\times8$&$1\times4$ & 15M& &$16\times16$& Total \\
3005 & 3004 & 3009 & 3001 & 3008 & 2431 & 32278  & 2780  & 6098 & \\
\hline
0.01\euro&0.01\euro&0.01\euro&0.01\euro&0.1\euro&0.01\euro&0.09\euro&0.01\euro&0.66\euro& \\
2$\times$&2$\times$&5$\times$&4$\times$&3$\times$&2$\times$&1$\times$&4$\times$&1$\times$&24 \\
\hline
0.02\euro&0.02\euro&0.05\euro&0.04\euro&0.03\euro&0.02\euro&0.09\euro&0.04\euro&0.66\euro&0.97\euro\\
\end{tabular}
\caption{Shopping list of the LEGO bricks. Buying components in different colors may be more expensive. Prices
are listed in  \euro , as of 30 November 2020, at \url{brickowl.com}.}
\label{tab:set}
} % \scriptsize
\end{center}
\end{table}

\begin{figure}[!htbp]
\begin{center}
    \fbox{\includegraphics[width=0.5\textwidth]{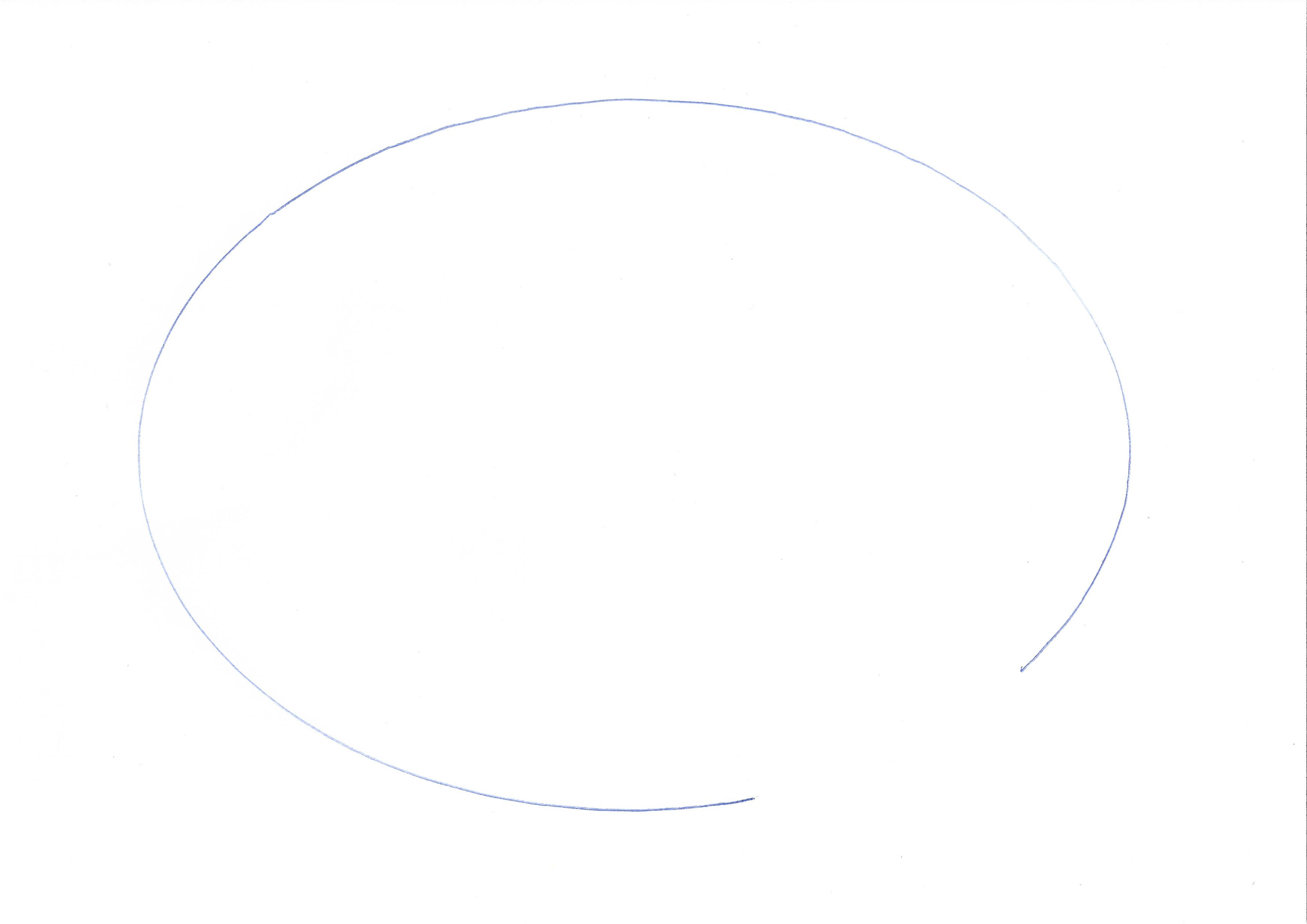}}
    \caption{An ellipse drawn with the ellipsograph}
    \label{ellipsograph4}
\end{center}
\end{figure}

\begin{figure}
\begin{center}
    \includegraphics[width=0.5\textwidth]{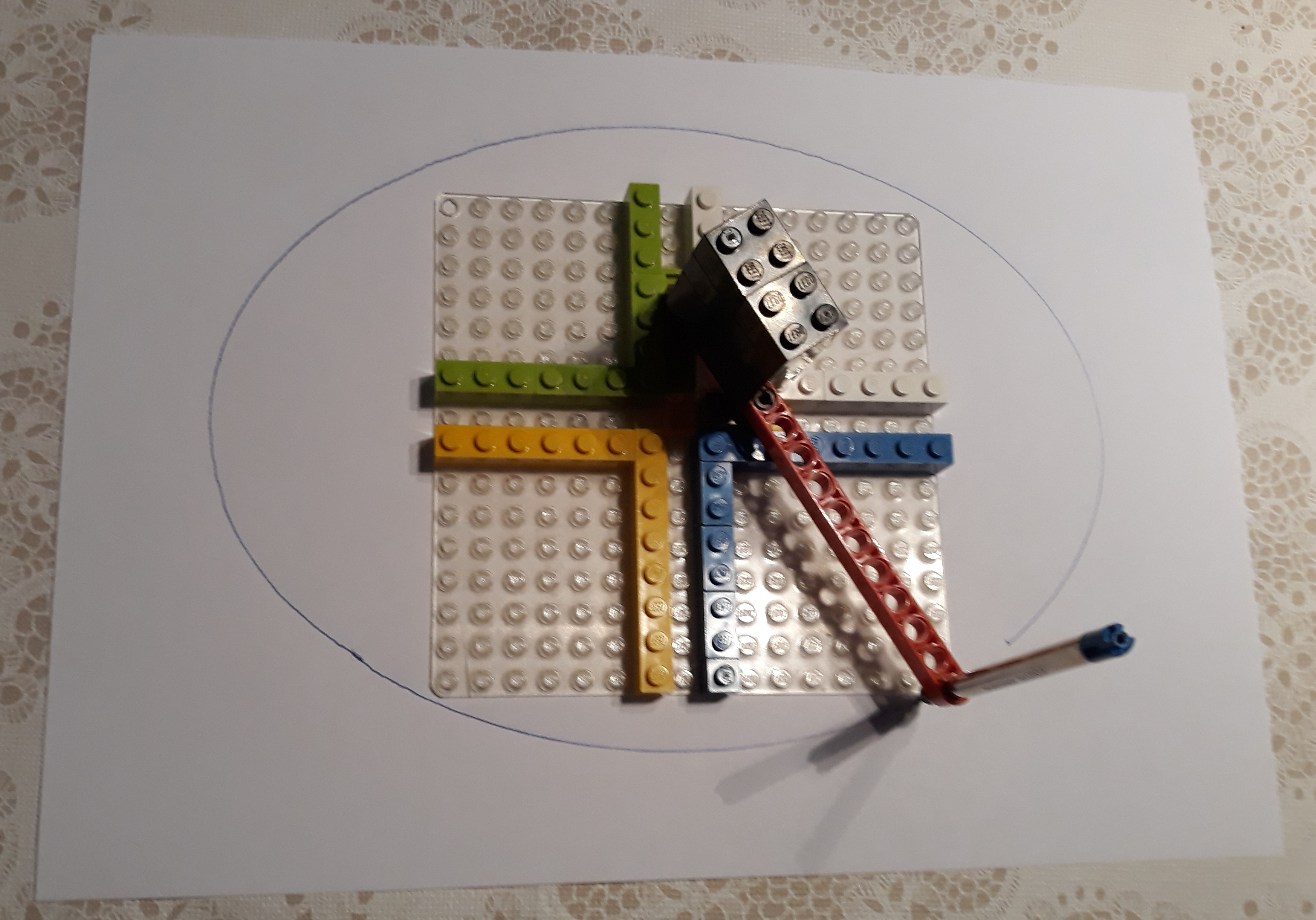}
    \caption{The ellipse and the ellipsograph}
    \label{ellipse-ellipsograph}
\end{center}
\end{figure}

Our contribution can draw almost 100\% of a quite big ellipse that fits on an A4 sheet of paper as it can be observed in Fig.~\ref{ellipsograph4} and \ref{ellipse-ellipsograph}. A very small part of the arc cannot be drawn because the shuttles get too close in that part of the movement. As future work we address solving this problem, too.

Finally, we point the reader to LEGO's official computer aided design tool, LEGO Digital Designer, available at \url{https://www.lego.com/en-us/ldd}. We used this program to visualize our concept in a digital way. The LXF file that contains our work can be downloaded at \url{https://matek.hu/zoltan/eg.lxf}. Also the building instructions can be checked out at \url{https://matek.hu/zoltan/Building%20Instructions%20%5Beg%5D.html}.
\begin{figure}[!htbp]
\centering
    \includegraphics[width=0.5\textwidth]{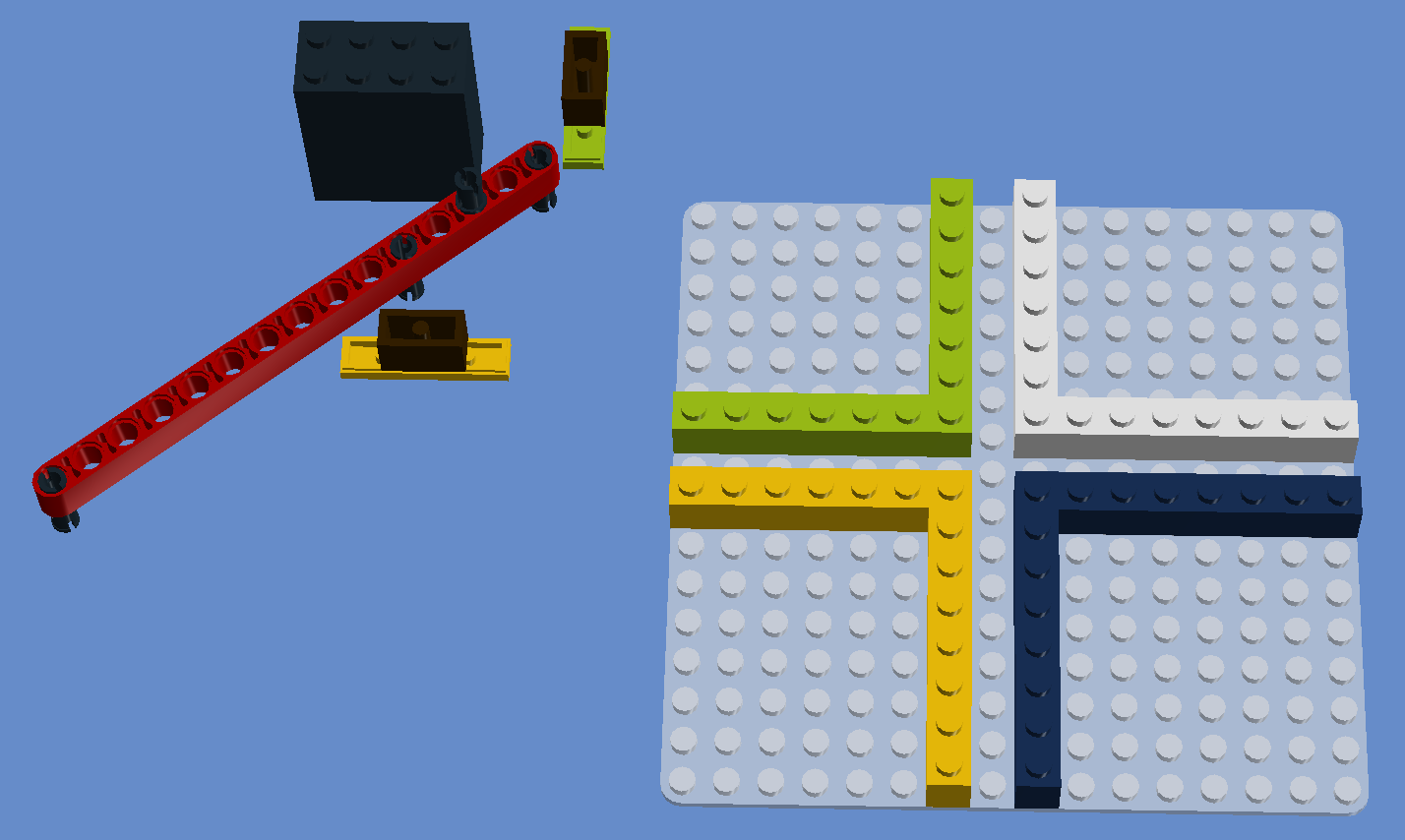}
    \caption{The ellipsograph sketched up with LDD}
    \label{ellipsograph5}
\end{figure}

\paragraph{Acknowledgements.}
Second author was partially supported by a grant MTM2017-88796-P from the
Spanish MINECO (Ministerio de Economia y Competitividad) and the ERDF
(European Regional Development Fund).

\bibliographystyle{plain} % We choose the "plain" reference style
\bibliography{ellipsograph} % Entries are in the "refs.bib" file

\end{document}